\documentclass[12pt, twoside, eqno]{article}
\usepackage{latexsym}
\usepackage{amssymb}
\usepackage{amsfonts}
\begin{document}

\noindent {\bf \large On ordered
$\Gamma$-semigroups ($\Gamma$-semigroups)}\bigskip

{\small \footnotetext{{\it Department of Mathematics, University of 
Athens, 15784
Panepistimiopolis, Athens, Greece}

email: nkehayop@math.uoa.gr}
\medskip

\noindent{\bf Niovi Kehayopulu}
\bigskip

\begin{abstract}
\noindent We add here some further characterizations to the 
characterizations of strongly regular ordered $\Gamma$-semigroups 
already considered in Hacettepe J. Math. 42 (2013), 559--567. Our 
results generalize the characterizations of strongly regular ordered 
semigroups given in the Theorem in Math. Japon. 48 (1998), 213--215, 
in case of ordered $\Gamma$-semigroups. The aim of writing this paper 
is not just to add a publication in $\Gamma$-semigroups but, mainly, 
as a continuation of the paper "On regular duo 
$po$-$\Gamma$-semigroups in Math. Slovaca 61 (2011), 871--884, to 
publish a paper which serves as an example to show what a 
$\Gamma$-semigroup is and give the right information about this 
structure.
\end{abstract}{\small


} }
\section {Introduction and prerequisites}
\noindent We have already seen in [2, 3] the methodology we use to 
pass from ordered semigroups (semigroups) to ordered 
$\Gamma$-semigroups ($\Gamma$-semigroups). Some of the results can be 
transferred just putting a "Gamma" in the appropriate place, while 
there are results for which the transfer is not easy. But anyway, we 
never work directly on ordered $\Gamma$-semigroups. If we want to get 
a result on an ordered $\Gamma$-semigroup, then we have to prove it 
first in an ordered semigroup and then we have to be careful to 
define the analogous concepts in case of the ordered 
$\Gamma$-semigroup (if they do not defined directly) and put the 
$"\Gamma"$ in the appropriate
place. In that sense, although a result on ordered $\Gamma$-semigroup 
generalizes its corresponding one of ordered semigroup, we can never 
say "we obtain and establish some important results in ordered 
$\Gamma$-semigroups extending and generalizing those for semigroups" 
(as we have seen in the bibliography) and this is because
we can do nothing on ordered $\Gamma$-semigroups if we do not examine 
it first for ordered semigroups. The same holds if we consider 
$\Gamma$-semigroups instead of ordered $\Gamma$-semigroups. In the 
present paper, the transfer was rather difficult, because it was not 
easy to give the definition of strongly regular ordered 
$\Gamma$-semigroups. An ordered $\Gamma$-semigroup $(M,\Gamma,\le)$ 
is called strongly regular if for every $a\in M$ there exist $x\in 
M$
and $\gamma, \mu\in \Gamma$ such that $a\le a\gamma x\mu a \mbox {
and } a\gamma x=x\gamma a=x\mu a=a\mu x$. This concept has been first 
introduced in Hacettepe J. Math. 42 (2013), 559--567, where this type 
of ordered $\Gamma$-semigroup has been characterized as an ordered 
$\Gamma$-semigroup $M$ which is both left and  right regular and the 
set $(M\Gamma a\Gamma M]$ is a strongly regular subsemigroup of $M$ 
for every $a\in M$. As a continuation of that result, we add here 
some further characterizations of this type of ordered 
$\Gamma$-semigroups, the main being the characterization of a 
strongly regular ordered $\Gamma$-semigroup as an ordered 
$\Gamma$-semigroup $M$ in which the $\cal N$-class $(a)_{\cal N}$ is 
a strongly regular subsemigroup of $M$ for every $a\in M$. So we 
prove that not only part of the Theorem in [5], but the whole Theorem 
can be transferred which is consistent  with what we already said in 
[2, 3]. The results of the present paper generalize the corresponding 
results on strongly regular ordered semigroups considered in [5].

A semigroup $S$ is called regular if for every $a\in S$ there exists 
$x\in S$ such that $a=axa$, that is, if $a\in aSa$ for every $a\in S$ 
or $A\subseteq ASA$ for every $A\subseteq S$. A semigroup $S$ is 
called left regular if for every $a\in S$ there exists $x\in S$ such 
that $a=xa^2$, that is $a\in Sa^2$ for every $a\in S$ or $A\subseteq 
SA^2$ for every $A\subseteq S$. It is called right regular if for 
every $a\in S$ there exists $x\in S$ such that $a=a^2x$, that is 
$a\in a^2S$ for every $a\in S$ or $A\subseteq A^2S$ for every 
$A\subseteq S$. A semigroup $S$ is called completely regular if for 
every $a\in S$ there exists $x\in S$ such that $a=axa$ and $ax=xa$ 
[6]. It has been proved in [6] that a semigroup is completely regular 
if and only if it is at the same time regular, left regular and right 
regular.

When we pass from semigroups to ordered semigroups, to completely 
regular semigroup correspond two concepts: The completely regular and 
the strongly regular ordered semigroups and this is the difference 
between semigroups and ordered semigroups.
For an ordered semigroup $S$ we denote by $(H]$ the subset of $S$ 
defined by $(H]:=\{t\in S \mid t\le h \mbox { for some } h\in H\}$. 
An ordered semigroup $(S,.,\le)$ is called regular if for every $a\in 
S$ there exists $x\in S$ such that $a\le axa$, that is, if $(a\in 
aSa]$ for every $a\in S$ or $A\subseteq (ASA]$ for every $A\subseteq 
S$. It is called left (resp. right) regular if for every $a\in S$ 
there exists $x\in S$ such that $a\le xa^2$ (resp. $a\le a^2 x$). An 
ordered semigroup $S$ is regular if and only if $a\in (aSa]$ for 
every $a\in S$ or $A\subseteq (ASA]$ for every $A\subseteq S$. It is 
left regular if and only if $a\in (Sa^2]$ for every $a\in S$ or 
$A\subseteq (SA^2]$ for every $A\subseteq S$. It is right regular if 
and only if $a\in (a^2S]$ for every $a\in S$ or $A\subseteq(A^2S]$ 
for every $A\subseteq S$. An ordered semigroup $S$ is called 
completely regular if it is regular, left regular and right regular. 
An ordered semigroup $S$ is called strongly regular if for every 
$a\in S$ there exists $x\in S$ such that $a\le axa$ and $ax=xa$. The 
strongly regular ordered semigroups are clearly completely regular 
but the converse statement does not hold in general. 
Characterizations of completely regular semigroups have been given in 
[6], characterizations of strongly regular ordered semigroups have 
been given in [5].

For two nonempty sets $M$ and $\Gamma$, define
$M\Gamma M$ as the set of all elements of the form $m_1\gamma
m_2$, where $m_1,m_2\in M$ and $\gamma \in \Gamma$. That is,

$M\Gamma M:=\{m_1 \gamma m_2 \mid m_1,m_2\in M, \gamma\in
 \Gamma\}$.

\noindent Let now $M$ and $\Gamma$ be two nonempty sets. The set
$M$ is called a {\it $\Gamma$-semigroup} if the following assertions 
are
satisfied:

(1) $M\Gamma M\subseteq M$.

(2) If $m_1,m_2,m_3,m_4\in M$, $\gamma_1,\gamma_2\in \Gamma$ such
that $m_1=m_3$, $\gamma_1=\gamma_2$

\hspace{0.6cm}and $m_2=m_4$, then $m_1\gamma_1m_2=m_3\gamma_2m_4$.

(3) $(m_1\gamma_1m_2)\gamma_2m_3=m_1\gamma_1(m_2\gamma_2m_3)$ for
all $m_1,m_2,m_3\in M$ and all
$\gamma_1,\gamma_2\in\Gamma$.\smallskip

\noindent This is the definition in [1] and it is a revised version 
of the definition of $\Gamma$-semigroups
given by Sen and Saha in [7], which allows us in an expression
of the form, say $A_1\Gamma A_2\Gamma,....., A_n\Gamma$ to put the
parentheses anywhere beginning with some $A_i$ and ending in some
$A_j$ $(i,j\in N=\{1,2,.....,n\})$ or in an expression of the form
$a_1\Gamma a_2\Gamma,....., a_n\Gamma$ or $a_1\gamma
a_2\gamma,....., a_n\gamma$ to put the parentheses anywhere
beginning with some $a_i$ and ending in some $a_j$ $(A_1,
A_2,....., A_n$ being subsets and $a_1,a_2,....., a_n$ elements of
$M$). Unless the uniqueness condition (widely used and still in use 
by some authors) in an expression of the form, say $a\gamma b\mu
c\xi d\rho e$ or $a\Gamma b\Gamma c\Gamma d\Gamma e$, it is not
known where to put the parentheses.

Here is an example of a $po$-$\Gamma$-semigroup $M$
which is easy to check [3] and shows exactly what a 
$\Gamma$-semigroup is. Other examples in which $M$ has order
3, 5 or 6 and $\Gamma$ order 2, one can find in [1--3]:
 Consider the two-elements set
$M:=\{a, b\}$, and let  $\Gamma=\{\gamma, \mu\}$ be the set of two
binary operations on $M$ defined in the tables
below:\begin{center} $\begin{array}{*{20}c}
   \gamma  &\vline &  a &\vline &  b  \\
\hline
   a &\vline &  a &\vline &  b  \\
\hline
   b &\vline &  b &\vline &  {a\,}  \\
\end{array}\,\,\,\,\,\,\,\,\,\,\,\,\,\,\,\,\,\,\,\,\,\,\,\,\,\,\begin{array}{*{20}c}
   \mu  &\vline &  a &\vline &  b  \\
\hline
   a &\vline &  b &\vline &  a  \\
\hline
   b &\vline &  a &\vline &  b  \\
\end{array}$
\end{center}
One can check that $(x\rho y)\omega z = x\rho (y\omega z)$ for all
$x,y,z\in M$ and all $\rho, \omega\in \Gamma$. So $M$ is a
$\Gamma$-semigroup.

An {\it ordered $\Gamma$-semigroup} (shortly {\it 
$po$-$\Gamma$-semigroup}) $M$, also denoted by\newline 
$(M,\Gamma,\le)$, is a $\Gamma$-semigroup $M$ endowed with an order 
relation $"\le"$ such that $a\le b$ implies $a\gamma c\le b\gamma c$ 
and $c\gamma a\le c\gamma b$ for all $c\in M$ and all 
$\gamma\in\Gamma$ [8]. An equivalence relation $\sigma$ on $M$ is 
called {\it congruence} if $(a,b)\in\sigma$ implies $(a\gamma 
c,b\gamma c)\in\sigma$ and $(c\gamma a,c\gamma b)\in\sigma$ for every 
$c\in M$ and every $\gamma\in\Gamma$. A congruence $\sigma$ on $M$ is 
called {\it semilattice congruence} if (1) $(a\gamma a,b\gamma 
a)\in\sigma$ for all $a,b\in M$ and all $\gamma\in\Gamma$ and (2) 
$(a,a\gamma a)\in\sigma$ for all $a\in M$ and all $\gamma\in\Gamma$.
A $po$-$\Gamma$-semigroup $M$ is called {\it regular} if $a\in 
(a\Gamma M\Gamma a]$ for every $a\in M$, equivalently if $A\subseteq 
(A\Gamma M\Gamma A]$ for every $A\subseteq M$. Keeping the already 
existing definition of left and right regular $\Gamma$ (or ordered 
$\Gamma$)-semigroups in the bibliography, we will call an ordered 
ordered $\Gamma$-semigroup $M$ {\it left} (resp. {\it right}) {\it 
regular} if $a\in (M\Gamma a\Gamma a]$ (resp. $a\in (a\Gamma a\Gamma 
M])$ for every $a\in M$. As in an ordered semigroup, we call a 
$po$-$\Gamma$-semigroup {\it completely regular} if it is at the same 
time regular, left regular and right regular. A nonempty subset $A$ 
of $M$ is called a {\it subsemigroup} of $M$ if $a,b\in A$ and 
$\gamma\in \Gamma$ implies $a\gamma b\in A$, that is if $A\Gamma 
A\subseteq A$.
A nonempty subset $A$ of $(M,\Gamma,\le)$ is called a {\it left} 
(resp. {\it right}) {\it ideal} of $M$ if (1) $M\Gamma A\subseteq A$ 
(resp. $A\Gamma M\subseteq A)$ and (2) if $a\in A$ and $M\ni b\le a$ 
implies $b\in A$. Clearly, the left ideals as well as the right 
ideals of $M$ are subsemigroups of $M$. A subsemigroup $F$ of 
$(M,\Gamma,\le)$ is called a {\it filter} of $M$ if (1) $a,b\in M$ 
and $\gamma\in\Gamma$ such that $a\gamma b\in F$ implies $a,b\in F$ 
and (2) if $a\in F$ and $M\ni b\ge a$ implies $b\in F$. For an 
element $a$ of $M$ we denote by $N(a)$ the filter of $M$ generated by 
$a$ and by $\cal N$ the relation on $M$ defined by ${\cal N}:=\{(a,b) 
\mid N(a)=N(b)\}$. Exactly as in ordered semigroups one can prove 
that the relation $\cal N$ is a semilattice congruence on $M$. A 
subset $T$ of an ordered $\Gamma$-semigroup $M$ is called {\it 
semiprime} if $A\subseteq M$ such that $A\Gamma A\subseteq T$ implies 
$A\subseteq T$, equivalently if $a\in M$ such that $a\Gamma 
a\subseteq T$ implies $a\in T$. For a $po$-$\Gamma$-semigroup $M$, we 
clearly have $M=(M]$, and for any subsets $A$, $B$ of $M$, we have 
$A\subseteq (A]=((A]]$, if $A\subseteq B$ then $(A]\subseteq (B]$, 
$(A]\Gamma (B]\subseteq (A\Gamma B]$ and $((A]\Gamma (B]]= ((A]\Gamma 
B]=(A\Gamma (B]]=(A\Gamma B]$.
\section{Main results}
\noindent{\bf Definition.} [4] A $po$-$\Gamma$-semigroup 
$(M,\Gamma,\le)$ is called
{\it strongly regular} if for every $a\in M$ there exist $x\in M$
and $\gamma, \mu\in \Gamma$ such that$$a\le a\gamma x\mu a \mbox {
and } a\gamma x=x\gamma a=x\mu a=a\mu x.$$A subsemigroup $T$ of 
$(M,\Gamma,\le)$ is called strongly regular if the set $T$ with the 
same $\Gamma$ and the order $"\le"$ of $M$ is strongly regular, that 
is, for every $a\in T$ there exist $y\in T$ and $\lambda, \rho\in 
\Gamma$ such that $a\le a\lambda y\xi a$ and $a\lambda y=y\lambda 
a=y\xi a=a\xi y$. We write it also as $(T,\Gamma,\le)$.\medskip

\noindent{\bf Theorem.} {\it Let M be an ordered
$\Gamma$-semigroup. The following are equivalent:
\begin{enumerate}
\item[$(1)$] M is strongly regular.
\item[$(2)$] For every $a\in M$, there exist $y\in M$ and 
$\gamma, \mu\in
\Gamma$ such that $a\le a\gamma y\mu a,\; y\le y\mu a\gamma 
y\mbox
{ and } a\gamma y=y\gamma a=y\mu a= a\mu y$.
\item[$(3)$] Every $\cal N$-class
if M is a strongly regular subsemigroup of M.
\item[$(4)$] The left and the right ideals of M are semiprime and 
for
every left ideal L and every right ideal R of M, the set 
$(L\Gamma
R]$ is a strongly regular subsemigroup of M.
\item[$(5)$] M is left regular, right regular, and the set 
$(M\Gamma
a\Gamma M]$ is a strongly regular subsemigroup of M for every
$a\in M$.
\item[$(6)$] For every $a\in M$ there exist $e_{a}$, 
$e^{'}_{a}\in M\Gamma
a\Gamma a\Gamma M$ and $\rho,\mu\in\Gamma$ such that $e_a \le
e_a\rho e^{'}_a$, $\,a\le e_a \mu a$, $\,a\le a\rho e^{'}_a$,
$(M\Gamma e_a\Gamma M]=(M\Gamma e^{'}_a\Gamma M]=(M\Gamma 
a\Gamma
M]$, and the set $(M\Gamma a\Gamma M]$ is a strongly regular 
subsemigroup
of $M$.
\item[$(7)$] For every $a\in M$ there exist $e_{a}$, 
$e^{'}_{a}\in M$ and $\rho,\mu\in\Gamma$ such that $\,a\le
e_a \mu a$, $\,a\le a\rho e^{'}_a$, and the set $(M\Gamma 
a\Gamma
M]$ is a strongly regular subsemigroup of $M$.
\item[$(8)$] For every $a\in M$, we have $a\in (M\Gamma a]\cap
(a\Gamma M]$, and $(M\Gamma a\Gamma M]$ is a strongly regular
subsemigroup of $M$.
\end{enumerate} }
\noindent{\bf Proof.} $(1)\Longrightarrow (2)$. For its proof we 
refer to [4]. For convenience, we sketch the proof: Let $a\in M$. 
Since $M$ is strongly regular, there exist $x\in M$ and $\gamma,
\mu\in \Gamma$ such that $a\le a\gamma x\mu a$ and $a\gamma
x=x\gamma a=x\mu a=a\mu x$. Then we have $$a\le a\gamma x\mu a\le
(a\gamma x\mu a)\gamma x\mu a=a\gamma (x\mu a\gamma x)\mu a.$$For the 
element $y:=x\mu a\gamma x$ of $M$, we have $a\le a\gamma y\mu
a$, $y\le y\mu a\gamma y$ and $a\gamma y=y\gamma a=y\mu a=a\mu 
y$.\smallskip

\noindent$(2)\Longrightarrow (3)$. Let $b\in M$. The class
$(b)_{\cal N}$ is a subsemigroup of $M$. Indeed: First of all, it
is a nonempty subset of $M$. Let $x,y\in (b)_{\cal N}$ and
$\gamma\in \Gamma$. Since $(x,b)\in\cal N$, $(b,y)\in\cal N$, and
$\cal N$ is a semilattice congruence on $M$, we have $(x\gamma y,
b\gamma y)\in\cal N$, $(b\gamma y, y\gamma y)\in\cal N$, $(y\gamma
y, y)\in\cal N$, then $(x\gamma y, y)\in\cal N$, and $x\gamma y
\in (y)_{\cal N}=(b)_{\cal N}$.\smallskip

\noindent$(b)_{\cal N}$ is strongly regular. In fact: Let $a\in
(b)_{\cal N}$. By (2), there exist $y\in M$ and $\gamma, \mu\in
\Gamma$ such that $a\le a\gamma y\mu a,\; y\le y\mu a\gamma y
\mbox { and } a\gamma y=y\gamma a=y\mu a= a\mu y$. On the other
hand, $y\in (b)_{\cal N}$. Indeed: Since $N(a)\ni a\le a\gamma
y\mu a$ and $N(a)$ is a filter of $M$, we have $a\gamma y\mu a\in 
N(a)$,
$y\in N(a)$, $N(y)\subseteq N(a)$. Since $N(y)\ni y\le y\mu
a\gamma y$ and $N(y)$ is a filter, we have $y\mu a\gamma y\in
N(y)$, $a\in N(y)$, $N(a)\subseteq N(y)$. Then $N(a)=N(y)$,
$(a,y)\in \cal N$, and $y\in (a)_{\cal N}=(b)_{\cal N}$.\smallskip

\noindent$(3)\Longrightarrow (4)$. Let $L$ be a left ideal of $M$
and $a\in M$ such that $a\Gamma a\subseteq L$. Then $a\in L$. In 
fact:
Since $a\in (a)_{\cal N}$ and $(a)_{\cal N}$ is strongly regular,
there exist $x\in (a)_{\cal N}$ and $\gamma, \mu\in\Gamma$ such
that $a\le a\gamma x\mu a$ and $a\gamma x=x\gamma a=x\mu a=a\mu
x$. Then we have
$$a\le (a\gamma x)\mu a=(x\gamma a)\mu a=x\gamma (a\mu
a)\in M\Gamma (a\Gamma a)\subseteq M\Gamma L\subseteq
L,$$and $a\in L$. If $R$ is a right ideal of $M$, $a\in M$ and 
$a\Gamma a\subseteq R$, then$$a\le a\gamma (x\mu a)=a\gamma (a\mu 
x)=(a\gamma a)\mu x\in (a\Gamma a)\Gamma M\subseteq R\Gamma 
M\subseteq R,$$so $a\in R$, and $R$ is also semiprime.

Let $L$ be a left ideal and $R$ a right ideal of $M$.
The (nonempty) set $(L\Gamma R]$ is a subsemigroup of $M$. In
fact: Let $a,b\in (L\Gamma R]$ and $\gamma\in\Gamma$. We have
$a\le y_1 \gamma_1 x_1$ and $b\le y_2 \gamma_2 x_2$ for some
$y_1,y_2\in L$, $\gamma_1,\gamma_2\in\Gamma$, $x_1,x_2\in R$. Then
$a\gamma b\le (y_1 \gamma_1 x_1)\gamma (y_2\gamma_2 x_2)$. Since
$(y_1 \gamma_1 x_1)\gamma y_2 \in(M\Gamma M)\Gamma L\subseteq
M\Gamma L\subseteq L$, we have $(y_1 \gamma_1 x_1\gamma y_2)\gamma
_2 x_2\in L\Gamma R$, so $a\gamma b\in (L\Gamma R]$.

\noindent Let now $a\in (L\Gamma R]$. Then there exist $x\in
(L\Gamma R]$ and $\gamma,\mu\in\Gamma$ such that $a\le a\gamma
x\mu a$ and $a\gamma x=x\gamma a=x\mu a=a\mu x$. In fact: Since
$a\in (a)_{\cal N}$ and $(a)_{\cal N}$ is strongly regular, there
exist $t\in (a)_{\cal N}$ and $\gamma,\mu\in\Gamma$ such that
$a\le a\gamma t\mu a$ and $a\gamma t=t\gamma a=t\mu a=a\mu t$.
Since $a\in (L\Gamma R]$, there exist $y\in L$, $\rho\in\Gamma$,
$z\in R$ such that $a\le y\rho z$. We have$$a\le a\gamma t\mu a\le
a\gamma t\mu (a\gamma t\mu a)=a\gamma (t\mu a\gamma t)\mu a.$$For
the element $x:=t\mu a\gamma t$, we have$$x=t\mu a\gamma t\le t\mu
(y\rho z)\gamma t=(t\mu y)\rho (z\gamma t).$$Since $t\mu y\in
M\Gamma L\subseteq L$ and $z\gamma t\in R\Gamma M\subseteq R$, we
have $(t\mu y)\rho (z\gamma t)\in L\Gamma R$, then $x\in L\Gamma
R$. Moreover, we have

$a\gamma x=x\gamma a$, that is, $a\gamma (t\mu a\gamma t)=(t\mu
a\gamma t)\gamma a$. Indeed,\begin{eqnarray*}a\gamma (t\mu a\gamma
t)&=&(a\gamma t)\mu (a\gamma t)=(t\mu a)\mu (t\gamma a)=t\mu (a\mu
t)\gamma a\\&=&t\mu (a\gamma t)\gamma a=(t\mu a\gamma t)\gamma
a.\end{eqnarray*}$x\mu a=a\mu x$, that is, $(t\mu a\gamma t)\mu
a=a\mu (t\mu a\gamma t)$. Indeed,\begin{eqnarray*}(t\mu a\gamma
t)\mu a&=&(t\mu a)\gamma (t\mu a)=(a\mu t)\gamma (a\gamma t)=a\mu
(t\gamma a)\gamma t\\&=&a\mu (t\mu a)\gamma t=a\mu (t\mu a\gamma
t).\end{eqnarray*}

$x\gamma a=x\mu a$, that is, $(t\mu a\gamma t)\gamma a=(t\mu
a\gamma t)\mu a$. Indeed,$$(t\mu a\gamma t)\gamma a=(t\mu a)\gamma
(t\gamma a)=(t\mu a)\gamma (t\mu a)=(t\mu a\gamma t)\mu a.$$
\noindent$(4)\Longrightarrow (5)$. Let $a\in M$. The set $(M\Gamma 
a\Gamma a]$ is a left ideal of $M$. This is because it is a nonempty 
subset of $M$ and we have\begin{eqnarray*}M\Gamma (M\Gamma a\Gamma 
a]&=&(M]\Gamma (M\Gamma a\Gamma a]\subseteq (M\Gamma M\Gamma a\Gamma 
a]=((M\Gamma M)\Gamma a\Gamma a]\\&\subseteq&(M\Gamma a\Gamma 
a]\end{eqnarray*}and $((M\Gamma a\Gamma a]]=(M\Gamma a\Gamma a]$.
Since $(M\Gamma
a\Gamma a]$ is a left ideal of $M$, by (4), it is semiprime. Since
$(a\Gamma a)\Gamma (a\Gamma a)\subseteq (M\Gamma a\Gamma a]$, we
have $a\Gamma a\subseteq (M\Gamma a\Gamma a]$, and $a\in (M\Gamma 
a\Gamma a]$, thus $M$ is left regular. Similarly the set $(a\Gamma 
a\Gamma M]$ is a right ideal of $M$ and $M$ is right regular. For the 
rest of the proof, we prove that $(M\Gamma
a\Gamma M]=((M\Gamma a]\Gamma (a\Gamma M]]$. Then, since $(M\Gamma
a]$ is a left ideal and $(a\Gamma M]$ a right ideal of $M$, by
(4), the set $((M\Gamma a]\Gamma (a\Gamma M]]$ is a strongly
regular subsemigroup of $M$, and so is $(M\Gamma a\Gamma M]$. We
have\begin{eqnarray*}M\Gamma a\Gamma M&\subseteq&M\Gamma (M\Gamma
a\Gamma a]\Gamma M=(M]\Gamma (M\Gamma a\Gamma a]\Gamma
(M]\\&\subseteq& ((M\Gamma M)\Gamma a\Gamma a\Gamma M]\subseteq
(M\Gamma a \Gamma a\Gamma M]\\&=&((M\Gamma a]\Gamma (a\Gamma
M]],\end{eqnarray*}then\begin{eqnarray*}(M\Gamma a\Gamma
M]&\subseteq&(((M\Gamma a]\Gamma (a\Gamma M]]]=((M\Gamma a]\Gamma
(a\Gamma M]]\\&=&((M\Gamma a)\Gamma (a\Gamma M)]\subseteq (M\Gamma
a\Gamma M],\end{eqnarray*}and so $(M\Gamma a\Gamma M]=((M\Gamma
a]\Gamma (a\Gamma M]]$.\smallskip

\noindent$(5)\Longrightarrow (6)$. Let $a\in M$. Since $M$ is left
regular, we have $a\in (M\Gamma a\Gamma a]$, since $M$ is right
regular, $a\in (a\Gamma a\Gamma M]$. Then there exist $x,y\in M$
and $\gamma, \mu,\rho,\xi\in\Gamma$ such that $a\le x\gamma a\mu
a$ and $a\le a\rho a\xi y$. Let
$e_a:=x\gamma a\rho a\xi y$ and $e^{'}_a:=x\gamma a\mu a\xi y$.
Then $e_a, e^{'}_a\in M\Gamma a\Gamma a\Gamma M$ and we have
$$a\le x\gamma a\mu a\le x\gamma (a\rho a\xi y)\mu
a=(x\gamma a\rho a\xi y)\mu a=e_a \mu a,$$
$$a\le a\rho a\xi y\le a\rho (x\gamma a\mu a)\xi y=a\rho (x\gamma
a\mu a\xi y)=a\rho e^{'}_a,$$
$$e_a=x\gamma a\rho a\xi y\le x\gamma
(a\rho a\xi y)\rho (x\gamma a\mu a)\xi y=e_a\rho
e^{'}_a.$$Moreover,\begin{eqnarray*} (M\Gamma e_a \Gamma
M]&=&(M\Gamma x\gamma a\rho a\xi y \Gamma M]\subseteq (M\Gamma
M\Gamma M\Gamma a\Gamma M\Gamma M]\\&\subseteq& (M\Gamma a\Gamma
M]\subseteq (M\Gamma (e_a\Gamma M]\Gamma M]=(M\Gamma e_a\Gamma
M\Gamma M]\\&\subseteq&(M\Gamma e_a\Gamma M],\end{eqnarray*}so
$(M\Gamma e_a\Gamma M]=(M\Gamma a\Gamma M]$, similarly $(M\Gamma
e^{'}_a\Gamma M]=(M\Gamma a\Gamma M]$. In addition, by (5), $(M\Gamma 
a\Gamma M]$ is a strongly regular subsemigroup of $M$.\smallskip

\noindent$(6)\Longrightarrow (7)$. This is obvious.\smallskip

\noindent$(7)\Longrightarrow (8)$. Let $a\in M$. By hypothesis,
there exist $e_{a}$, $e^{'}_{a}\in M$ and
$\rho,\mu\in\Gamma$ such that $\,a\le e_a \mu a$,
$\,a\le a\rho e^{'}_a$, and the set $(M\Gamma a\Gamma M]$ is a 
strongly
regular subsemigroup of $M$. Since $a\le e_a\mu a$, we have $a\in
(M\Gamma a]$. Since $a\le a\rho e^{'}_a\in a\Gamma M$, we have
$a\in (a\Gamma M]$. Then $a\in (M\Gamma a]\cap (a\Gamma
M]$ and $(M\Gamma a\Gamma M]$ is a strongly regular subsemigroup of 
$M$, so condition (8) is satisfied. \smallskip

\noindent$(8)\Longrightarrow (1)$. For its proof we refer to [4]. 
$\hfill\Box$\medskip

\noindent{\bf Corollary.} {\it Let M be an ordered 
$\Gamma$-semigroup. The following are equivalent:
\begin{enumerate}
\item[$(1)$] M is strongly regular.
\item[$(2)$] If $a\in M$, then for the subset $e_a:=M\Gamma 
a\Gamma a\Gamma M$ of $M$, we have

    $e_a\subseteq (e_a\Gamma e_a]$, $\;a\in (e_a\Gamma a]$, $\; 
a\in (a\Gamma e_a]$, $(M\Gamma e_a \Gamma M]=(M\Gamma a\Gamma 
M]$,

    and the set $(M\Gamma a\Gamma M]$ is a strongly regular 
subsemigroup of M.
\item[$(3)$] If $a\in M$, then there exists a subset $e_a$ of $M$ 
such that

    $\;a\in (e_a\Gamma a]$, $\; a\in (a\Gamma e_a]$,

     and the set $(M\Gamma a\Gamma M]$ is a strongly regular 
subsemigroup of M.
\end{enumerate} }

\noindent{\bf Proof.} $(1)\Longrightarrow (2)$. Let $a\in M$ and 
$e_a:=M\Gamma a\Gamma a\Gamma M$. Since $M$ is strongly regular, by 
the Theorem $(1)\Rightarrow (5)$, $M$ is left regular, right regular, 
and the set $(M\Gamma a\Gamma M]$ is a strongly regular subsemigroup 
of $M$. Since $M$ is left regular and right regular, we have $a\in 
(M\Gamma a\Gamma a]$ and $a\in (a\Gamma a\Gamma M]$. Then we 
have\begin{eqnarray*}e_a&=&M\Gamma a\Gamma a\Gamma M\subseteq M\Gamma 
(a\Gamma a\Gamma M]\Gamma (M\Gamma a\Gamma a]\Gamma M\\&=&(M]\Gamma 
(a\Gamma a\Gamma M]\Gamma (M\Gamma a\Gamma a]\Gamma 
(M]\\&\subseteq&((M\Gamma a\Gamma a\Gamma M)\Gamma (M\Gamma a\Gamma 
a\Gamma M)]\\&=&(e_a\Gamma e_a],\end{eqnarray*}
\begin{eqnarray*}a\in (M\Gamma a\Gamma a]&\subseteq&(M\Gamma (a\Gamma 
a\Gamma M]\Gamma a]=
((M\Gamma (a\Gamma a\Gamma M)\Gamma a]\\&=&(e_a\Gamma 
a],\end{eqnarray*}
\begin{eqnarray*}a\in (a\Gamma a\Gamma M]&\subseteq&(a\Gamma (M\Gamma 
a\Gamma a]\Gamma M]=
(a\Gamma (M\Gamma a\Gamma a\Gamma M)]\\&=&(a\Gamma 
e_a],\end{eqnarray*}
\begin{eqnarray*}(M\Gamma e_a\Gamma M]&=&(M\Gamma(M\Gamma a\Gamma 
a\Gamma M)\Gamma M]\subseteq (M\Gamma a\Gamma M]\\&\subseteq&(M\Gamma 
(M\Gamma a\Gamma a]\Gamma M]=(M\Gamma (M\Gamma a\Gamma a)\Gamma 
M]\\&\subseteq& (M\Gamma (M\Gamma a\Gamma (a\Gamma a\Gamma M]\Gamma 
M]\\&=&(M\Gamma (M\Gamma a\Gamma a\Gamma a\Gamma M)\Gamma 
M]\\&\subseteq& (M\Gamma (M\Gamma a\Gamma a\Gamma M)\Gamma M]\\&=&
(M\Gamma e_a\Gamma M],
\end{eqnarray*}thus we have $(M\Gamma e_a\Gamma M]=(M\Gamma a\Gamma 
M]$.\\$(2)\Longrightarrow (3)$. This is obvious.\smallskip

\noindent$(3)\Longrightarrow (1)$. Let $a\in M$. By hypothesis, there 
exists a subset $e_a$ of $M$ such that $a\in (e_a\Gamma a]$, $a\in 
(a\Gamma e_a]$, and the set $(M\Gamma a\Gamma M]$ is a strongly 
regular subsemigroup of $M$. Since $a\in (e_a\Gamma a]\subseteq 
(M\Gamma a]$ and $a\in (a\Gamma e_a]\subseteq (a\Gamma M]$, we have
$a\in (M\Gamma a]\cap (a\Gamma M]$. By the Theorem $(8)\Rightarrow 
(1)$, $M$ is strongly regular. $\hfill\Box$\medskip

\noindent{\bf Remark.} Here are some information we get about ordered 
semigroups: By the implication $(1)\Rightarrow (2)$ of the above 
Corollary (or by the implication $(1)\Rightarrow (2)$ of the Theorem 
in [5]), we have the following: If $S$ is a strongly regular ordered 
semigroup, $a\in S$ and $e_a:=Sa^2S$, then we have $e_a\subseteq 
({e_a}^2]$, $a\in (e_a a]$, $a\in (ae_a]$, $(Se_a S]=(SaS]$, and the 
set $(SaS]$ is a strongly regular subsemigroup of $S$. By the 
implication $(3)\Rightarrow (1)$ of the Corollary (or by the 
implication $(8)\Rightarrow (1)$ of the Theorem in [5]), we have the 
following: Let $S$ be an ordered semigroup. Suppose that for every 
$a\in S$ there exists a subset $e_a$ of $S$ such that $a\in (e_a a]$, 
$a\in (ae_a]$, and the set $(SaS]$ is a strongly regular subsemigroup 
of $S$. Then $S$ is strongly regular.{\small

\end{document}